\documentclass[centertags,12pt]{amsart}
\usepackage{latexsym}
\usepackage{amssymb}
\usepackage{graphics}

\numberwithin{equation}{section}
\makeatletter
\renewcommand{\subsection}{\@startsection
{subsection}{2}{0mm}{\baselineskip}{-0.25cm}
{\normalfont\normalsize\bf}}
\makeatother

\newtheorem{theorem}{Theorem}[section]
\newtheorem{proposition}[theorem]{Proposition}
\newtheorem{corollary}[theorem]{Corollary}
\newtheorem{lemma}[theorem]{Lemma}
{\theoremstyle{definition}

\newtheorem{example}[theorem]{Example}
}
\theoremstyle{remark}
\newtheorem{remark}[theorem]{Remark}

\def\1{\mathbf 1}

\def\F{\mathbf F}

\def\N{\mathbf N}

\def\cC{\mathcal C}

\def\cF{\mathcal F}

\def\cL{\mathcal L}
\def\cP{\mathcal P}
\def\cX{\mathcal X}

\def\cV{\mathcal V}

\def\supp{{\rm Supp}}
\def\div{{\rm div}}
\def\Div{{\rm Div}}
\def\dim{{\rm dim}}
\def\deg{{\rm deg}}

\begin{document}

\author[C.~Munuera]{Carlos Munuera}
\author[F.~Torres]{Fernando Torres}

\thanks{The first author was supported by the Grant VA020-02 from the ``Junta de Castilla y Le\'on". This paper was written while the second author was visiting The University of Valladolid (Dpto. de Algebra, Geometr\'{\i }a y Topolog\'{\i }a; Fac. Ciencias) supported by the Grant SB2000-0225 from the ``Secretaria de Estado de Educaci\'on y Universidades del Ministerio de Educaci\'on, Cultura y Deportes de Espa\~na"}

\title[A bound on the trellis state complexity of AG codes]{A Goppa-like bound on the Trellis state complexity of algebraic
geometric codes}

\address{Dept. of Applied Mathematics, University of Valladolid (ETS Arquitectura), Avda. Salamanca SN, 47014 Valladolid,
Castilla, Spain} \email{cmunuera@modulor.arq.uva.es}

\address{IMECC-UNICAMP, Cx. P. 6065, Campinas, 13083-970-SP, Brazil}
\email{ftorres@ime.unicamp.br} 

\address{Current Address: Departamento de Algebra, Geometr\'{\i }a y Topolog\'{\i }a, Facultad de Ciencias - Universidad
de Valladolid, c/ Prado de la Magdalena s/n 47005, Valladolid (Spain)} \email{ftorres@agt.uva.es}

    \begin{abstract} For a linear code $\cC$ of length $n$ and dimension $k$, Wolf noticed that the trellis state
complexity $s(\cC)$ of $\cC$ is upper bounded by $w(\cC):=\min(k,n-k)$. In this paper we point out some new lower bounds
for $s(\cC)$. In particular, if $\cC$ is an Algebraic Geometric code, then $s(\cC)\geq w(\cC)-(g-a)$, where $g$ is the
genus of the underlying curve and $a$ is the abundance of the code.
    \end{abstract}

\maketitle

    \section{Introduction}\label{s1}

Let $n$ be a positive integer. A {\em trellis} of {\em depth} $n$ is an edge-labeled directed graph $T=(V,E)$ with vertex
set $V$ and edge set $E$, where $V$ is the disjoint union of $(n+1)$ sub-sets $V_0,V_1,\ldots,V_n$, such that every edge
in $E$ that begins at $V_i$ ends at a vertex in $V_{i+1}$, and every vertex in $V$ belongs to at least one path from a
vertex in $V_0$ to a vertex in $V_n$. In this paper we only consider trellises with $V_0$ and $V_n$ having just one
element.

Although the trellises are graph-theoretic objects, they play an important role in Coding Theory. Firstly, their use was
restricted to convolutional coding. Later they were used also for block codes, mainly for the purpose of soft-decision
decoding with the Viterbi algorithm. History and the state of the art of applications of trellises to Coding Theory can be
seen in Forney \cite{forney} or Vardy's survey \cite{vardy}.

Now to each path from $V_0$ to $V_n$ of a trellis $T$ of depth $n$ we can associate an ordered $n$-tuple over a label
alphabet, say $\F$ a finite field. We say that the trellis $T$ {\em represents} a block code $\cC$ of length $n$ over
$\F$, if the set of all such $n$-tuples read from $T$ is exactly the set of codewords of $\cC$. There might exist more than one
non-isomorphic trellises that represents the same code. A way to measure the {\em complexity} of the trellis $T$ that
represents $\cC$ is by means of the {\em state complexity} of $T$, defined by    
    $$
s_T(\cC):=\max(s_0(T),s_1(T),\ldots,s_n(T))\, ,
    $$
where $s_i(T)={\rm log}_q|V_i|$, being $q=|\F|$ and $V_0,V_1,\ldots,V_n,$ the underlying partition of the vertex set of
$T$. Moreover, if the code $\cC$ is linear, and once its immersion in $\F^n$ (i.e., the order of coordinates of $\cC$) is
fixed, there exists a unique (up to a graph isomorphism) trellis $T_\cC$ such that for each $i=0,1,\ldots,n$, and any
trellis $T$ that represents $\cC$ it holds that $s_i(T_\cC)\leq s_i(T)$. The trellis $T_\cC$ is called the {\em minimal
trellis} of $\cC$. However, as was noticed by Forney \cite{forney} (see also \cite[Ex. 5.1]{vardy}), the state of
complexity of $T_\cC$ depends on the immersion of $\cC$ in $\F^n$. Thus one is lead to consider the so-called {\em trellis
state complexity} of $\cC$, denoted by $s(\cC)$, as being the minimal state of complexity of $T_\cC$ under any immersion
of $\cC$ in $\F^n$. We point out that the role of the state of complexity in the study of a linear code is comparable to
the role that plays its length, its dimension and its minimum distance; cf. Muder \cite{muder}, Forney \cite{forney}.

There are several bounds on $s(\cC)$ available in the literature. Here we just mention the Wolf bound, as it was first
noticed by him in \cite{wolf}, namely    
    $$
s(\cC)\leq w(\cC):=\min(k,n-k)\, ,
    $$
where $k$ is the dimension of $\cC$. For further upper and lower bounds on the trellis state complexity of linear codes 
see \cite[Sects. 5.2, 5.3]{vardy}.

The study of the trellis state complexity of some types of codes (including BCH, RS, and RM codes) has been carried out by
several authors; see \cite{berger-beery1, berger-beery2, blackmore-norton2, kasami-all, vardy-beery}. The case of
Algebraic Geometric codes (or simply, AG codes) was treated by Blackmore and Norton \cite{blackmore-norton1}, and by Shany
and Be'ery \cite{shany-beery}. Their results mainly concern on AG codes on the Hermitian curve; nevertheless, these
authors noticed that the trellis state complexity of an AG code
$\cC(D,G)$ on a curve of genus $g$ attains the Wolf bound
provided that either $\deg(G)<\lfloor \deg(D)/2\rfloor$, or $\deg(G)>\lceil\deg(D)/2\rceil+2g-2$. Otherwise, the trellis
state complexity is lower bounded by $\lceil\deg(D)/2\rceil-g-1$; see Proposition \ref{prop3.1} and Corollary
\ref{cor3.3} here.

The purpose of this paper is to point out some new lower bounds for the trellis state complexity of a linear code. Propositions \ref{prop2.1} and \ref{prop2.2} show lower bounds on $s(\cC)$ for any linear code $\cC$. These bounds are related to either the minimum distance or to generalized Hamming weights of $\cC$, and particular cases of them have been already noticed by Forney \cite[Sect. IIIC]{forney} and by Vardy and Be'ery \cite{vardy-beery}. Their applications to AG codes led to an interesting lower bound on $s(\cC)$ which has to do with the gonality sequence of the curve used to construct the code (see Proposition \ref{prop3.2}). As a consequence we obtain the aforementioned results of Blackmore-Norton and Shany-Be'ery as well as the main result of this paper, namely Theorem \ref{thm3.1} which for an AG code $\cC$ says that
     $$
s(\cC)\geq w(\cC)-(g-a)\, ,
     $$
where $g$ is the genus of the underlying curve and $a$ is the abundance of the code. This bound is analogous to the Goppa bound on the minimum distance $d$ of $\cC$, namely (see Corollary \ref{cor3.1})
     $$
d\geq S(\cC)-(g-a)\, ,
     $$
where $S(\cC):=n-k+1$ is the Singleton bound on the minimum distance of $\cC$. We illustrate our results via Examples \ref{ex3.1} and \ref{ex3.2}.

     \section{Lower bounds on the trellis state complexity of a linear
code}\label{s2}

Throughout this section we let $\cC$ be a $[n,k,d]$ linear code over a finite field $\F$. Let $\cC^\perp$ be its dual with
minimum distance $d^\perp$. For a fixed immersion of $\cC$ (and the corresponding of $\cC^\perp$) in $\F^n$, we let $T$
and $T^\perp$ be the minimal trellises of $\cC$ and $\cC^\perp$ respectively.

    \begin{lemma}{\rm (\cite[Thm. 4.20]{vardy})}\label{lemma2.1}\quad $s_T(\cC)=s_{T^\perp}(\cC^\perp).$
    \end{lemma}

From this lemma and Forney's construction below, the Wolf bound $s(\cC)\leq w(\cC):=\min(k,n-k)$ follows.

Now, according to Forney's construction of the minimal trellis $T$ of $\cC$ (see \cite{forney}), we can assume that the
sub-sets $V_0,V_1,\ldots,V_n$ of the underlying partition of the vertex set of $T$ are given by  
   $$
V_i=\cC/\cP_i\oplus \cF_i\, ,
   $$
where $\cP_i=\cP_i(\cC)$ and $\cF_i=\cF_i(\cC)$ are respectively the $i$-th past and the $i$-th future sub-codes of $\cC$, namely $\cP_0=\cF_n=0$, $\cP_n=\cF_0=\cC$ and for $i=1,\ldots,n-1$:
    \begin{align*}
\cP_i & =  \{(c_1,\ldots,c_i): (c_1,\ldots,c_i,0,\ldots,0)\in\cC\}\, ,\\         
\cF_i & =  \{(c_{i+1},\ldots,c_n): (0,\ldots,0,c_{i+1},\ldots,c_n)
\in\cC\}\, .                                   
    \end{align*}
It follows that $|V_i|$ is a power of $|\F|$ so that
    $$
s_i=s_i(T)=k-\dim_\F(\cP_i)-\dim_\F(\cF_i)\qquad\text{for each $i=0,1,\ldots,n$}\, .
    $$  
    \begin{remark}\label{rem2.0} By the definition of $d$, $\cP_i=0$ for $i=1,\ldots,d-1$ and $\cF_i=0$ for $i=n-d+1,\ldots,n-1$. Therefore
    $$
s_T(\cC)=\max(s_{d-1},\ldots,s_{n-d+1})
    $$
whenever $2d\leq n+3$; otherwise, $s_T(\cC)=k$.
    \end{remark}
    \begin{proposition}\label{prop2.1} For a non-negative integer $t$, suppose that either $2d\geq n+2-t,$ or $2d^\perp\geq n+2-t.$ Then $s(\cC)\geq w(\cC)-t.$
    \end{proposition}
    \begin{proof} As we already remarked, $\cP_{d-1}=(0)$. Now suppose that $2d\geq n+2-t$. If $\cF_{d-1}=(0)$, then $s_{d-1}(T)=k$ and hence $s_T(\cC)\geq k$. If $\cF_{d-1}\neq (0)$, then $\cF_{d-1}$ is a linear code of length at most $(n-d+1)$ whose minimum distance is at least $d$. Thus, by applying the Singleton bound to $\cF_{d-1}$, we obtain
   $$
\dim_\F(\cF_{d-1})\leq (n-d+1)-d+1=n-2d+2\leq t\, ,
   $$
and hence $s_T(\cC)\geq s_{d-1}(T)\geq k-t$. Now, if $2d^\perp\geq n+2-t$, from the above proof we obtain $s_{T^\perp}(\cC^\perp)\geq n-k-t$ so that $s_T(\cC)\geq w(\cC)-t\ (*)$ by Lemma \ref{lemma2.1}. Thus the result follows since the bound $(*)$ does not depend of the immersion of $\cC$ in $\F^n$.
   \end{proof}
   \begin{remark}\label{rem2.1} If $\cC$ is an MDS code, then it always hold that $\max(2d,2d^\perp)\geq n+2$. Hence, as a
corollary of the above proposition, we get the well known result that MDS codes reach the Wolf bound.
   \end{remark}
Next we use the idea of the above proof in order to derive another new lower bound on $s(\cC)$. This  bound has to do
with the generalized Hamming weight hierarchy $(d_i(\cC):i=1,\ldots,k)$ of $\cC$, where
    $$
d_i(\cC):=\min\{\#\supp(\cV): \text{$\cV$ is a subcode of $\cC$ of dimension at least $i$}\} 
    $$
is the so-called $i$-th generalized Hamming weight of $\cC$. Notice that $d_1(\cC)=d$.
    \begin{proposition}\label{prop2.2} For an integer $i$ with $1\leq i\leq k,$ suppose that either $d_i(\cC)\geq n+2-d,$
or $d_i(\cC^\perp)\geq n+2-d^\perp$. Then $s(\cC)\geq w(\cC)-i+1.$
    \end{proposition}
    \begin{proof} Let $d_i(\cC)\geq n+2-d$. Since the length of $\cF_{d-1}$ is at most $n-d+1$, then 
$\#\supp(\cF_{d-1})<d_i(\cC)$ by the hypothesis on $d_i(\cC)$. Thus $\dim_\F(\cF_{d-1})<i$ by the definition of
$d_i(\cC)$,
and since $\cP_{d-1}=(0)$, we have that $s_T(\cC)\geq k-i+1$. Similarly, if $d_i(\cC^\perp)\geq n+2-d^\perp$, then
$s_{T^\perp}(\cC^\perp)\geq n-k-i+1$ and the result follows.
    \end{proof}

    \section{Lower bounds on the trellis state complexity of an AG
code}\label{s3}

In this section we state and prove some new lower bounds on the trellis state complexity of an AG code which will be related to the gonality sequence of the curve used to construct the code. Standard references on Algebraic Geometry and Algebraic Function Fields and Codes are respectively Hartshorne's book \cite{hartshorne} and Stichtenoth's book \cite{sti}.

An AG code (also called a Geometric Goppa code) is constructed from the following data:
    \begin{itemize}
\item A projective, geometrically irreducible, non-singular algebraic curve $\cX$ defined over a finite field $\F$ (or simply, a curve over $\F$);
\item Two $\F$-divisors on $\cX$: $D=P_1+\ldots+P_n$ and $G$, such that their support are disjoint, $P_i\neq P_j$ for $i\neq j$, and the point  $P_i$ is an $\F$-rational point of $\cX$ for all $i=1,\ldots,n$.
    \end{itemize}
Then the AG code $\cC_\cX(D,G)$ is the image in $\F^n$ of the $\F$-linear map
    $$
ev: \cL(G)\to \F^n\, ,\qquad f\mapsto (f(P_1),\ldots,f(P_n)\, ,
    $$
where $\cL(G)=\{f\in\F(\cX)^*:\div(f)+G\succeq 0\}\cup\{0\}$. The dimension $k$ and the minimum distance $d$ of $\cC_\cX(D,G)$ can be estimated by means of the Riemann-Roch theorem as they satisfy:
    \begin{itemize}
\item $k=\ell(G)-\ell(G-D)$;
\item $d\geq n-\deg(G)$ (the Goppa bound).
    \end{itemize}
We will always assume that the divisor $(G-D)$ is special, otherwise  $\cC_\cX(D,G)=\F^n$. The dimension over $\F$ of the
kernel of the map $ev$, $a:=\ell(G-D)$, is called the {\em abundance} of $\cC_\cX(D,G)$. The dual code of $\cC_\cX(D,G)$
is also an AG code, namely $\cC_\cX(D,D-G+K)$ where $K$ is an $\F$-canonical divisor on $\cX$ obtained as the divisor of
a differential form having simple poles and residue $1$ at each point in $\supp(D)$.
     \begin{remark}\label{rem3.1} Let $\cC:=\cC_\cX(D,G)$. If $\deg(G)\geq n$, then the Goppa bound on the minimum distance of $\cC$ does not give any
information. In this case, Corollary \ref{cor3.0} below gives an improvement on the Goppa bound by means of the gonality sequence of $\cX$ provided that $a>0$; as a matter of fact, this corollary  gives an improvement on a previous result due to Pellikaan \cite[Thm. 2.11]{pellikaan0}. On the other hand, if $n>\deg(G)$ (so that $a=0$) and $g$ is the genus of $\cX$,
then by the Riemann-Roch theorem the Goppa bound implies
     $$
d\geq S(\cC)-g\, ,
     $$
where $S(\cC)=n-k+1$ is the Singleton bound on $\cC$; see Corollary \ref{cor3.1} for the case $a>0$.
    \end{remark}
Next we recall some basic facts on the gonality sequence of curves. For the curve $\cX$ over $\F$ and a positive integer
$i$, let 
    $$
\gamma_i=\gamma_i(\cX,\F):=\min\{\deg(A): \text{$A\in\Div_\F(\cX)$ and 
                               $\ell(A)\geq i$}\}\, .
    $$
The sequence $GS(\cX)=GS(\cX,\F):=(\gamma_i:i\in\N)$ is called the {\em gonality sequence} of $\cX$ over $\F$; 
cf. \cite{yang-kumar-sti}. Notice
that $\gamma_1=0$ and that $\gamma_2$ is the usual gonality of $\cX$ over $\F$. Some properties of $GS(\cX)$ are stated below.    
     \begin{lemma}{\rm (\cite[Prop. 11]{yang-kumar-sti})}\label{lemma3.1} Let $g$ be the genus of $\cX$ and suppose that
$\cX(\F)\neq \emptyset.$ Then
     \begin{enumerate}
\item[\rm(1)] The sequence $GS(\cX)$ is strictly increasing$;$
\item[\rm(2)] $\gamma_{i}=g+i-1$ for $i\geq g+1;$
\item[\rm(3)] $\gamma_g=2g-2;$
\item[\rm(4)] $\gamma_i\geq 2i-2$ for $i=1,\ldots,g-1.$
     \end{enumerate}
     \end{lemma}
     \begin{lemma}{\rm (\cite[Cor.~2.4]{pellikaan})}\label{lemma3.2} If $\cX$ is a plane curve of degree $r$ and
$\cX(\F)\neq\emptyset,$ then $GS(\cX)$ is the strictly increasing sequence obtained from the semigroup generated by $r-1$
and $r.$
     \end{lemma}
The Goppa bound is generalized to higher Hamming weights as follows.
    \begin{lemma}{\rm (\cite[Cor.~2]{munuera})}\label{lemma3.3} The $i$-th generalized Hamming weight $d_i$ of the
AG code $\cC_\cX(D,G)$ satisfies $d_i\geq n-\deg(G)+\gamma_{a+i},$ where $a$ is the abundance of the code and
$\gamma_{a+i}$ is the $(a+i)$-th element of the gonality sequence of $\cX.$
    \end{lemma}
     \begin{proof} (Sketch) We have that $d_i\leq t$ if and only if there exist $(n-t)$ pairwise different points in
$\supp(D)$, say
$P_{i_1},\ldots, P_{i_{n-t}}$, such that $\ell(G-P_{i_1}-\ldots-P_{i_{n-t}})\geq a+i$. Thus
    $$
d_i=\min\{n-\deg(D'): \text{$0\preceq D'\preceq D$ and $\ell(G-D')\geq
a+i$}\}\, ,
    $$
and the result follows.
     \end{proof}
    \begin{corollary}{\rm (Improved Goppa bound)}\label{cor3.0} For any AG code $\cC_\cX(D,G)$ of abundance $a$, 
   $$
d\geq n-\deg(G)+\gamma_{a+1}\, .
   $$
    \end{corollary}
As we already noticed in Remark \ref{rem3.1}, the following result is well-known for $a=0$.
     \begin{corollary}\label{cor3.1} For any AG code $\cC_\cX(D,G)$ over a curve $\cX$ of genus $g$ and abundance $a,$ 
    $$
d\geq (n-k+1)-(g-a)\, .
    $$ 
     \end{corollary}
     \begin{proof} By the Riemman-Roch theorem, the dimension $k$ of the code satisfies $k\geq \deg(G)+1-g-a$. From this
inequality and Corollary \ref{cor3.0} we have
    $$
d\geq (n-k+1)-g+\gamma_{a+1}-a\, ,
    $$
and the results follows as $\gamma_{a+1}\geq 2a$ (see Lemma \ref{lemma3.1}) once we check that $a\leq g\ (*)$. In fact, we
can assume $a\geq 1$, and $(*)$ follows by applying Clifford's theorem to the divisor $(G-D)$ since by assumption it is
special so that $\deg(G-D)\leq 2g-2$.
    \end{proof}
Next we deal with lower bounds on the trellis complexity of AG codes. Throughout we let 
    $$
   \cC:=\cC_\cX(D,G)
    $$ 
be an AG code defined by the divisors $D=P_1+\ldots+P_n$ and $G$ on a curve $\cX$ over $\F$ of genus $g$. We let $k$ be
the dimension, $d$ the minimum distance, and $a=\ell(G-D)$ the abundance of $\cC$. We let $(\gamma_i:i\in\N)$ denote the
gonality sequence of $\cX$ over $\F$.

Our first two results are applications of Propositions \ref{prop2.1} and \ref{prop2.2} to the case of AG codes. 
    \begin{proposition}{\rm (\cite[Prop. 1]{shany-beery})}\label{prop3.1} We have $s(\cC)=w(\cC)$ provided that either
     \begin{enumerate}
\item[\rm(1)] $\deg(G)<\lfloor n/2\rfloor+\gamma_{a+1},$ or
\item[\rm(2)] $\deg(G)>\lceil n/2\rceil+2g-2-\gamma_{a+1}.$
     \end{enumerate}
     \end{proposition}
     \begin{proof} If (1) holds, from the improved Goppa bound (Corollary \ref{cor3.0}) we have that $d\geq n-\lfloor
n/2\rfloor+1\geq (n+2)/2$ and the result follows from Proposition \ref{prop2.1} with $t=0$. If (2) holds, we apply the
above argument to the dual code of $\cC$ and the proof is complete.
     \end{proof}
     \begin{proposition}\label{prop3.2} Let $i$ be a positive integer. Then $s(\cC)\geq w(\cC)-i+1$ provided that either
     \begin{enumerate}
\item[\rm(1)] $\gamma_{a+i}\geq 2\deg(G)-n-\gamma_{a+1}+2,$ or
\item[\rm(2)] $\gamma_{a+i}\geq n+2(2g-2)-2\deg(G)-\gamma_{a+1}+2.$
     \end{enumerate}
     \end{proposition}
    \begin{proof} If $i\geq k+1$, the result is clear; otherwise, it is a consequence of Proposition \ref{prop2.2}, Lemma
\ref{lemma3.3}, and the Improved Goppa bound Corollary \ref{cor3.0}.
    \end{proof}
Now we can state the main result of this paper.
     \begin{theorem}\label{thm3.1}\quad $s(\cC)\geq w(\cC)-(g-a).$
     \end{theorem}
     \begin{proof} It is enough to show that one of the hypotheses in Proposition \ref{prop3.2} is fulfilled with
$i_0=g+1-a$. Recall that $a\leq g$ (see the proof of Corollary \ref{cor3.1}) so that $i_0\geq 1$. By Lemma
\ref{lemma3.1}, $\gamma_{g+1}=2g$. Now 
     $$
2g\geq 2\deg(G)-n-\gamma_{a+1}+2\ \Leftrightarrow\ 2\deg(G)\leq n+2g-2+\gamma_{a+1}\, , \qquad\text{and}
     $$
     $$
2g\geq n+2(2g-2)-2\deg(G)-\gamma_{a+1}\ \Leftrightarrow\ 
2\deg(G)\geq n+2g-2-\gamma_{a+1}\, .
     $$
Therefore in order that one of the hypothesis in Proposition \ref{prop3.2} with $i_0=g+1-a$ holds true it is enough that 
     \begin{equation}\label{eq3.1}
n+2g-2-\gamma_{a+1}\leq n+2g-2+\gamma_{a+1}\, ,
     \end{equation}
and the result follows.
     \end{proof}
     \begin{remark}\label{rem3.2} If $a\leq g/2$, we can apply the above proof with $i_0=g+1-2a$ to show that
$s(\cC)\geq w(\cC)-(g-2a)$. In fact, in this case condition (\ref{eq3.1}) becomes
     $$
n+2g-2+2a-\gamma_{a+1}\leq n+2g-2-2a+\gamma_{a+1}\, ,
     $$
which is true as $\gamma_{a+1}\geq 2a$ by Lemma \ref{lemma3.1}.
     \end{remark}
     \begin{remark}\label{rem3.3} The bound obtained in Theorem \ref{thm3.1} is sharp in some cases as we shall see in Example \ref{ex3.1} below, where AG codes on rational and elliptic curves are considered.
     \end{remark}
     \begin{corollary}\label{cor3.2} Suppose that
    \begin{equation}\label{eq3.2}
\lfloor n/2\rfloor +\gamma_{a+1}\leq \deg(G)\leq \lceil n/2\rceil+2g-2-\gamma_{a+1}\, .
    \end{equation}
Then 
    $$
s(\cC)\geq w(\cC)-\min(\deg(G)+1-a-\lfloor (n+\gamma_{a+1})/2\rfloor, 2g-1-\deg(G)-a+\lceil (n-\gamma_{a+1})/2\rceil)\, . 
    $$
     \end{corollary}
     \begin{proof} Let $\alpha:=\deg(G)+1-a-\lfloor (n+\gamma_{a+1})/2\rfloor$ and $\beta:=2g-1-\deg(G)-a+\lceil
(n-\gamma_{a+1})/2\rceil$. We show that $s(\cC)\geq w(\cC)-\alpha\ (*)$. If $g-a\leq \alpha$, then $(*)$ follows from
Theorem \ref{thm3.1}. Let $g-a\geq \alpha+1$ and set $i_0:=\alpha+1$. We have that $i_0\geq 1$ by the inequality in the
left side of (\ref{eq3.2}), and  $\gamma_{a+i_0}\geq 2(a+i_0-1)$ by Lemma \ref{lemma3.1}. Then an
straightforward computation shows that hypothesis (1) in Proposition \ref{prop3.2} holds true and so $(*)$ follows from
that proposition. Similarly, we have that $s(\cC)\geq w(\cC)-\beta$ and the result follows.
     \end{proof}
The next result was obtained in \cite{blackmore-norton2} and
\cite[Prop. 1]{shany-beery} whenever $\lfloor
n/2\rfloor >2g-2$ and $a=0$.
     \begin{corollary}{\rm (Clifford bound, \cite[Prop. 1]{shany-beery})}\label{cor3.3} If (\ref{eq3.2}) holds, then
   $$
s(\cC)\geq \lfloor (n+\gamma_{a+1})/2\rfloor -g\, .
   $$
     \end{corollary}
     \begin{proof} Let $i_0$ be as in the above proof. Then from that proof and that of Proposition \ref{prop2.2} we have
that $s(\cC)\geq k-i_0+1$. Since $k=\ell(G)-a\geq \deg(G)+1-g-a$, the proof follows.
     \end{proof}
Next we point out two more lower bounds on $s(\cC)$ that follows from
the formulae below. Let us recall that 
the $i$-th past $\cP_i$ and the $i$-th future $\cF_i$  
sub-codes of $\cC$ in Forney's construction of the minimal trellis $T_\cC$
that represents the code $\cC$ are given respectively by (see
\cite{blackmore-norton1, shany-beery})
     \begin{align*}
\cP_i & =  \cC(D-P_{i+1}-\ldots-P_n,G-P_{i+1}-\ldots-P_n)\, ,\\         
\cF_i & =  \cC(D-P_1-\ldots-P_i,G-P_1-\ldots-P_i)\, ;                                   
    \end{align*}
in particular, for each $i=0,\ldots,n$,
    \begin{equation}\label{eq3.3}
s_i(T)=\ell(G)+a-\ell(A_i)-\ell(B_i)\, ,
    \end{equation}
where $A_i:=G-P_1-\ldots-P_i$ and $B_i:=G-P_{i+1}-\ldots-P_n$.
     \begin{proposition}\label{prop3.3}\quad $s(\cC)\geq k+2a-\ell(2G-D)-1.$
     \end{proposition}
     \begin{proof} If there exists $i\in\{0,1,\ldots,n\}$ such that 
$\ell(A_i)\geq 1$ and $\ell(B_i)\geq 1$, then $\ell(A_i)+\ell(B_i)\leq
\ell(A_i+B_i)+1$ (see e.g. \cite[Lemma
5.5]{hartshorne}) and the result follows from (\ref{eq3.3}). Suppose now that for all $i=0,1,\ldots,n$, either $\ell(A_i)=0$ or $\ell(B_i)=0$. Then $a=0$ and $\min\{\ell(A_i)+\ell(B_i):i=0,1,\ldots,n\}\leq 1$. Hence from (\ref{eq3.3}), $s_{T_\cC}(\cC)\geq k-1\geq k-\ell(2G-D)-1$ and the result follows.
      \end{proof}
     \begin{proposition}\label{prop3.4} Let $i\in\{0,1,\ldots,n\}$ and $j\in \N$ such that $\deg(G)-\gamma_j<\min(i,n-i).$ Then 
    $$
s(\cC)\geq k-2(j-1-a)\, .
    $$
    \end{proposition}
    \begin{proof} By hypothesis $\deg(A_i)<\gamma_j$ and
$\deg(B_i)<\gamma_j$. Then by the definition of $\gamma_j$, $\ell(A_i)\leq
j-1$ and $\ell(B_i)\leq j-1$, and the result follows from (\ref{eq3.3}).
     \end{proof}
     \begin{remark}\label{rem3.4} If $n$ is even, the hypothesis in the proposition above is equivalent to $\deg(G)-\gamma_j<n/2$.
     \end{remark}
Finally, we consider the case of selfdual AG codes.
    \begin{proposition}\label{prop3.5} Let $i\in\{0,1,\ldots,n\}$. For a
selfdual AG code $\cC_\cX(D,G)$ on a curve $\cX$ of genus $g$, 
   $$
s_i(T_\cC)=n-i-2(\ell(A_i)-a)\, .
   $$
    \end{proposition}
    \begin{proof} Since the code is selfdual, there is a canonical divisor
$K$ on $\cX$ such that $2G=D+K$. Thus
$\ell(B_i)=\ell(K-A_i)$ and the conclusion follows from the Riemann-Roch
theorem and (\ref{eq3.3}).
    \end{proof}
    \begin{corollary}\label{cor3.4} Let $\cC=\cC_\cX(D,G)$ be a selfdual
AG code on an elliptic curve $\cX.$ Then
    $$
s_i(T_\cC)=\begin{cases}
i        & \text{if $0\leq i<n/2$,} \\
n/2-2\ell(A_{n/2})    & \text{if $i=n/2$,} \\
n-i    & \text{if $n/2<i\leq n$.} 
    \end{cases}
    $$
In particular$,$ $s_{T_\cC}(\cC)=w(\cC)-1$ if and only if $G\sim P_1+\ldots+P_{n/2};$ otherwise $s_{T_\cC}(\cC)=w(\cC).$
     \end{corollary}
     \begin{proof} We have $\deg(G)=n/2$ so that $a=0$ and $w(\cC)=n/2$. Thus $s_i(T_\cC)=n-2\ell(A_i)-i$ for
$i=0,1,\ldots,n$ by
Proposition \ref{prop3.5}. If $i<n/2$, $\ell(A_i)=n/2-i$ and so $s_i(T_\cC)=i$; if $i>n/2$, $\ell(A_i)=0$ and so
$s_i(T_\cC)=n-i$. Finally if $i=n/2$, $\ell(A_{n/2})\in\{0,1\}$. Since $\ell(A_{n/2})=1$ if and only if $G\sim P_1+\ldots+P_{n/2}$, the result follows.
     \end{proof}
     \begin{remark}\label{rem3.5} Observe that the above result on $s_{T_\cC}(\cC)$ depends on the immersion of the code in
$\F^n$. The corresponding result for the trellis state complexity of a selfdual AG code $\cC=\cC_\cX(D,G)$ on an elliptic
curve $\cX$ is the
following: We have $s(\cC)=w(\cC)-1$ if and only if there exist $n/2$ pairwise different points
$P_{i_1},\ldots,P_{i_{n/2}}\in\supp(D)$ such that $G\sim P_{i_1}+\ldots+P_{i_{n/2}}$; otherwise $s(\cC)=w(\cC)$.
    \end{remark}
To end this paper we shall include some examples showing the performance of the bounds we have obtained. 
    \begin{example}\label{ex3.1} In this example we show that the Goppa bound on the trellis state complexity (Theorem \ref{thm3.1}) is sharp in some cases. More precisely, we consider AG codes on rational or elliptic curves. 

{\em Rational codes.} If $\cX$ is rational, then $\cC=\cC_\cX(D,G)$ is an MDS code; moreover $a=0$. Hence, according to 
Theorem \ref{thm3.1}, it holds that $s(\cC)=w(\cC)$.
 
{\em One-point elliptic codes.} Let $\cX$ be an elliptic curve over $\F$. Let $\{P_1,\dots,P_n, O\}$ be the set of $\F$-rational points of $\cX$, where $O$ is the neutral element of the group $(\cX(\F),\oplus)$. For $1< m< n-1$, let us consider the one-point  elliptic code $\cC(m)$ associated to the divisors $D=P_1+\ldots+P_n$ and $G=mO$. The current complexity measures of ${\mathcal C}(m)$ have been computed in \cite{eliptico}. In order to collect these results, we first need some notation. A point $P\in {\mathcal X}(\F)$ is called {\em idempotent} if $P\oplus P=O$. Let $\delta$ denote the number of idempotent points in ${\mathcal X}(\F)$. It is easily seen that $\delta\in\{ 1,2,4\}$. Set $\ell_1(m):=\ell(mO-P_{n-m+1}-\ldots -P_n)$ and  $\ell_2(m):=\ell(mO-P_{1}-\cdots -P_{m})$. Note that $0\leq \ell_1(m),\ell_2(m)\leq 1$. 

Let $T$ be a minimal trellis that represents ${\mathcal C}(m)$.
According to Wolf bound and Goppa bound (Theorem \ref{thm3.1}), its state complexity $s_T(\cC(m))$ satisfies $w(\cC(m))\geq s_T(\cC(m))\geq w(\cC(m))-1$. The exact values of $s_T(\cC(m))$ are as follows:

(1) If $n>2m$, then
    $$
s_T(\cC(m))=\begin{cases}
k          &  \text{if $n>2m+1$ or $\ell_1(m)+\ell_2(m)<2$,}\\
k-1        &  \text{if $n=2m+1$ and $\ell_1(m)+\ell_2(m)=2$.}
     \end{cases}
    $$
If either $n>2m+1$, or $n=2m+1$ and $\delta=2$, then there exists an ordering on $\supp(D)$ such that $\ell_1(m)+\ell_2(m)=2$. If $n=2m+1$ and $\delta=4$, then there is no ordering on $\supp(D)$ such that $\ell_1(m)+\ell_2(m)=2$.
    
(2) If $n=2m$ then
    $$
s_T(\cC(m))=\begin{cases}
k          &  \text{if $\ell_1(m)+\ell_2(m)=0$,}\\
k-1        &  \text{if $\ell_1(m)+\ell_2(m)\neq 0$.}
    \end{cases}
   $$
In this case, there always exists an ordering on $\supp(D)$ such that $\ell_1(m)+\ell_2(m)=1$.

(3) If $n<2m$ then
   $$
s_T(\cC(m))=\begin{cases}
n-k          &   \text{if $n<2m-1$ or $\ell_1(m)+\ell_2(m)<2$,}\\
n-k-1        &   \text{if $n=2m-1$ and $\ell_1(m)+\ell_2(m)=2$.}
    \end{cases}
   $$
There exists an ordering such that $\ell_1(m)+\ell_2(m)=2$ excepts in the case that $n=2m-1$ and $\delta=4$.
    \end{example}
    \begin{example}\label{ex3.2} Let $\cC=\cC_\cX(D,G)$ be an AG code on a curve $\cX$ of genus $g$ over $\F$ of cardinality $q$. Suppose that $\deg(G)$ satisfies inequalities (\ref{eq3.2}). (Otherwise, $s(\cC)=w(\cC)$ by Proposition \ref{prop3.1}.) Suppose in addition that
    \begin{equation}\label{eq3.4}
n=\deg(D)>2(2g-2)+1\, .
    \end{equation}
This hypothesis is satisfied if for instance $\cC$ is a one-point Goppa code, $n=|\cX(\F)|-1$ and the underlying curve $\cX$ being one the following:
    \begin{itemize}
\item $\F$-maximal; i.e., such that $|\cX(\F)|=q+1+2g\sqrt{q}$ (the maximum number of $\F$-rational points according to the Hasse-Weil bound);
\item the Suzuki curve, namely the one admitting the plane model $y^q-y=x^{q_0}(x^q-x)$, where $q_0=2^s, s\geq 1$, $q=2q_0^2$. Here $|\cX(\F)|=q^2+1$ and $g=q_0(q-1)$ (see \cite{hansen-sti}).
     \end{itemize}
A first consequence of (\ref{eq3.2}) and (\ref{eq3.4}) is that $a=0$ so that $k=\dim_\F(\cC)=\deg(G)+1-g$. In particular,
    $$
w(\cC)=\begin{cases}
\deg(G)+1-g & \text{if $\lfloor n/2\rfloor\leq \deg(G)\leq \lfloor n/2\rfloor +g-1$,}\\
n-\deg(G)+g-1 & \text{if $\lfloor n/2\rfloor +g\leq\deg(G)\leq\lceil n/2\rceil+2g-2$.}
      \end{cases}
    $$
Thus Corollaries \ref{cor3.3} and \ref{cor3.4} imply the same result, namely 
    $$
s(\cC)\geq \lfloor n/2\rfloor -g\, .
    $$
We know present two numerical examples.

(1) Let $\cX$ be the Hermitian curve $Y^5Z+YZ^5=X^6$ over $\F$ of cardinality $q=25$. It is $\F$-maximal of genus $g=10$; i.e. $|\cX(\F)|=126$. We let $D=P_1+\ldots+P_n$ with $n=125$, $G=mQ$ with $Q=(0:1:0)$ and $m$ satisfying (\ref{eq3.2}), i.e., $62\leq m\leq 81$. However by duality (cf. Lemma \ref{lemma2.1}) we can assume
   $$
62\leq m\leq 71\, .
   $$
Set $\cC(m):=\cC_\cX(D,mQ)$. In the following table we list lower bound on $s(\cC(m))$ obtained from our results. The column ``Wolf" means the value $w(\cC(m))$ which is computed above; the column ``Clifford" means $\lfloor n/2\rfloor -g=52$ (Corollary \ref{cor3.4}); the column ``Prop. \ref{prop3.2}" means the bound obtained from Proposition \ref{prop3.2}(1). Finally the columns ``Prop. \ref{prop3.3}" and ``Prop. \ref{prop3.4}" show the lower bounds obtained respectively from Propositions \ref{prop3.3} and \ref{prop3.4}.  Here we use Lemma \ref{lemma3.2} to compute the gonality sequence of $\cX$ (namely, $GS(\cX)=\{0,5,6,10,11,12,15,16,17,18,20,\to\}$), and use the facts that $\div(x^{25}-x)=D-125 Q$ and that the Weierstrass semigroup at $Q$ is generated by $5$ and $6$ (see e.g. \cite[Sect. III]{yang-kumar-sti}) in order to compute $\ell(2G-D)$ as required by Proposition \ref{prop3.3}. The best lower bound obtained in this way is written in bold face.
   \begin{center}
\begin{tabular}{||c||c|c|c|c|c|c|} 
\hline
$m$ & Wolf  & Clifford & Prop. \ref{prop3.2} & Prop. \ref{prop3.3} & Prop. \ref{prop3.4} \\ 
\hline

62  & 53    &  {\bf 52} & {\bf 52}      & {\bf 52}    & 51 \\

63  & 54    &      52   & {\bf 53}      & 52          & 52 \\

64  & 55    &      52   & {\bf 54}      & 53          & 53 \\

65  & 56    &      52   &   53          & 53          & {\bf 54} \\

66  & 57    &      52   &   54          & 53          & {\bf 55} \\

67  & 58    &      52   & {\bf 54}      & {\bf 54}    & {\bf 54} \\

68  & 59    &      52   & {\bf 53}      & {\bf 53}    & {\bf 53} \\

69  & 60    &      52   & {\bf 54}      & 53          & {\bf 54} \\

70  & 61    &      52   & 53            & 53          & {\bf 55} \\

71  & 62    &      52   & 52            & 52          & {\bf 56} \\ 
\hline
    \end{tabular}
    \end{center} 
(2) Now we let $\cX$ be the Suzuki curve over $\F$ of cardinality $q=2q_0^2$, $q_0=2$. Thus $|\cX(\F)|=65$ and $g=14$. Let $\cC(m)$ be the one-point Goppa code on $\cX$ defined by $D=P_1+\ldots+P_{64}$, $G=mQ$ with $Q$ the unique point over $x=\infty$ and $32\leq m\leq 45$. Here we cannot apply Lemma \ref{lemma3.2} to compute the gonality sequence of $\cX$ because it is defined by a singular plane curve. However, as $\div(x^8-x)=D-64 Q$, and the Weierstrass semigroup at $Q$ is generated by $8, 10, 13$ and $14$ (cf. \cite{hansen-sti}) we can compute the bound in Proposition \ref{prop3.3} and obtain the following table.
 \begin{center}
\begin{tabular}{||c||c|c|c|c|c|c|} 
\hline
$m$ & Wolf  & Clifford & Prop. \ref{prop3.3} \\ 
\hline

32  & 19    &  {\bf 18}      & 17   \\

33  & 20    &  {\bf 18}      & {\bf 18}   \\

34  & 21    &  18      & {\bf 19}   \\

35  & 22    &  18      & {\bf 20}   \\

36  & 23    &  18      & {\bf 20}   \\

37  & 24    &  18      & {\bf 20}   \\

38  & 25    &  18      & {\bf 20}   \\

39  & 26    &  18      & {\bf 20}   \\

40  & 27    &  18      & {\bf 20}   \\

41  & 28    &  18      & {\bf 20}   \\ 

42  & 29    &  18      & {\bf 20}   \\

43  & 30    &  18      & {\bf 20}   \\

44  & 31    &  18      & {\bf 20}   \\

45  & 32    &  18      & {\bf 20}   \\
\hline
    \end{tabular}
    \end{center} 
    \end{example}
    
\end{document}